\documentclass[12pt]{article}
\usepackage{amsmath}
\usepackage{amssymb}

\begin{document}

\title{Generalized complex structures on Kodaira surfaces}
\author{Vasile Brinzanescu\thanks{Corresponding author.} \footnote{Institute of Mathematics ``Simion Stoilow'' of the Romanian Academy, P.O. Box 1-764, RO-70700, Bucharest, Romania. Vasile.Brinzanescu@imar.ro} \  and\  Oana Adela Turcu\footnote{oannat@yahoo.com}}


\date{}
\maketitle

\begin{abstract}

We compute the deformations in the sense of generalized complex structures of the standard classical complex structure on a primary Kodaira surface and we prove that the obtained family of deformations is a smooth locally complete family depending on four complex parameters. This family is the same as the extended deformations (in the sense of Kontsevich and Barannikov) in  degree two, obtained by Poon using differential Gerstenhaber algebras.

\end{abstract}

\section*{1. Introduction }

In the paper [7], Nigel Hitchin introduced new geometrical structures which unify many classical structures  and sheds light on mirror symmetry and many of the supersymmetric equations arising from string and $M$-theory and supergravity. He defines a generalized complex structure to be a complex structure, not on the  tangent bundle $T$ of a manifold $N$, but on $T \oplus T^*$, unifying in this way the complex geometry and the  symplectic geometry in some sense. Indeed, a complex  or symplectic structure on a manifold determines a maximal isotropic subbundle of the bundle $T\oplus T^*$. 
On the smooth sections of the bundle $T \oplus T^*$ there is a natural bracket, called Courant bracket, and Hitchin defines a generalized complex structure as an almost complex structure $J$ on $T \oplus T^* \; (J^2 = -1)$, whose $i$-eigenbundle $L \subset T \oplus T^*$ is Courant involutive. This new geometrical structure is, in some sense, the complex analogue of the Dirac structure introduced by Courant and Weinstein   [4], [5], in order to unify Poisson geometry with symplectic geometry.  The study of generalized complex structures was continued by Gualtieri (see [6]).

In this paper we start with the (classical) complex structure on a primary Kodaira surface and we compute, by using properties of the Lie algebroids, the family of deformations of this complex structure in the sense of generalized complex structures (see [6]). By solving the generalized Maurer-Cartan equation we get the main result of the paper, which shows that the obtained family of deformations is a smooth locally complete family depending on four complex parmeters. This family is the same as the extended deformations (in the sense of Kontsevich [9] and Barannikov - Kontsevich [1]) in degree two, obtained by Poon in [12], by  using differential Gerstenhaber algebras. 

In particular, we get the family (depending of two complex parameters) of the deformations of  (classical) complex structures on a Kodaira surface, obtained by Borcea in [2], as well as examples of  generalized complex structures of complex type, which are not (classical) complex structures on  a (primary) Kodaira surface. 

\vspace{0.5cm}

{\it Acknowledgment} The first named author expresses his gratitude to the Max-Plank-Institut f\" ur Mathematik Bonn for its hospitality and stimulating atmosphere; part of this paper was prepared during his stay at the Institute. Partially supported by CNCSIS contract 1189 / 2009-2011.

\vspace{0.5cm}

\section*{2. Kodaira surfaces}

\vspace{0.5cm}

A primary Kodaira surface is a complex analytic fibre bundle of elliptic curves over an elliptic curve 
(see [8] or [3]). In fact, a Kodaira surface is an elliptic complex surface, 
with trivial canonical bundle, of the form $N = \mathbb{C}^2/G$, where $\mathbb{C}^2$ denotes the space of two complex variables $(z, w)$ and $G$ is  a properly discontinuous
non-abelian group of affine transformations, generated by four elements $g_i$, $i = 1, 2, 3,4$, satisfying the following relations:
$$
g_i g_j g_i^{-1} g_j^{-1} = \mbox{id},
$$
for all couples $(i,j)$, $i < j$, $(i, j) \neq (3, 4)$ and
$$
g_3g_4g_3^{-1}g_4^{-1} = g_2^m, 
$$
for some positive integer $m$. 

The universal covering of $N$ is complex analitically isomorphic to $\mathbb{C}^2$, the space of two complex variables $(z, w)$, such that
$g_j$ regarded as covering transformations ($G$ is the fundamental group of $N$), take the form:
$$
g_j(z, w) = (z + \alpha_j, w + \bar \alpha_j \cdot z + \beta_j), \quad j = 1,2,3,4,
$$
where $\alpha_1 = \alpha_2 = 0$ and $\bar \alpha_3 \alpha_4 - \bar\alpha_4\alpha_3 = m\beta_2$.

We shall identify $\mathbb{C}^2$ with $\mathbb{R}^4$, the space of four real variables $(x, y, u, v)$ by $z = x + iy$, $w = u + iv$. Now,  considering the group $A$ of all real-affine transformations of $\mathbb{C}^2 = \mathbb{R}^4$, commuting with any transformation of the form
$$
g(z, w) = (z + \alpha, w + \bar \alpha z + \beta), 
$$
one can identify $A$ to $\mathbb{C}^2$ endowed with the following multiplication: 
$$
(z, w) \star (\alpha, \beta) = (z + \alpha, w + \bar\alpha z + \beta).
$$

Thus, the complex structure of $\mathbb{C}^2$ determines a right-invariant complex structure on $A$ and the Kodaira surface $N$ appears as
the compact quotient of $A$ by some discrete subgroup $G$, generated by 
$$
g_j(0, 0) = (\alpha_j, \beta_j), \quad j = 1,2,3,4;
$$
(see, for example [2] or [12]). 

From the point of view  of differential (or even real-analytic) structure, a Kodaira surface is a  parallelizable manifold, i.e. the tangent bundle $T_N$ is globally generated by invariant vector fields $\{X, Y, U, V\}$; the only non-zero Poisson bracket is 
$[X, Y] = U$ (see [12] or [2]). The complex structure endomorphism $J$ is acting on $T_N$ by 
$$
JX = Y, \quad JY = -X, \quad JU = V, \quad JV = -U;
$$
(see [12]). 

\vspace{0.5cm}

\section*{3. Generalized complex structures on manifolds and deformations}

\vspace{0.5cm}

A generalized complex structure is a notion introduced  by Hitchin [7] and developed by Gualtieri [6]. Generalized complex geometry contains complex and symplectic geometry as extremal special cases. 

A generalized complex structure on a smooth manifold $M$ is defined to be a complex structure $J\; (J^2 = -1)$, not on the tangent bundle $T_M$, but on the sum $T \oplus T^*$ of the tangent and cotangent bundles, which is required to  be orthogonal with respect to the natural inner product on sections $X + \sigma, Y + \tau \in \mathcal{C}^\infty(T_M \oplus T_M^*)$ defined by 
$$
\langle X + \sigma, Y +\tau \rangle = \frac{1}{2} (\sigma (Y) + \tau(X)).
$$

This is only possible if $\dim_{\mathbb{R}} M = 2n$, which we suppose. 

In addition, the $(+ i)$-eigenbundle
$$
L \subset (T_M \oplus T_M^*) \otimes \mathbb{C}
$$
of $J$ is required to be involutive with respect to the Courant bracket, a skew bracket operation on smooth sections of $T_M \oplus T_M^*$ defined by
$$
[X + \sigma, Y + \tau] = [X, Y] + \mathcal{L}_X \tau - \mathcal{L}_Y\sigma - \frac{1}{2} d(i_X\tau - i_Y\tau),
$$
where $\mathcal{L}_X$ and $i_X$ denote the Lie derivative and interior product operations on forms. 

Since $J$ is orthogonal with respect to $\langle \cdot, \cdot\rangle$, the $(+i)$-eigenbundle $L$ is a maximal isotropic subbundle of $(T_M \oplus T_M^* )\otimes \mathbb{C}$ of real index zero (i.e. $L \cap \bar L = \{0\}$). 

In fact, a generalized complex structure on $M$ is completely determined by a maximal isotropic subbundle $L \subset (T_M \oplus T_M^*) \otimes \mathbb{C}$ of real index zero, which is  Courant involutive (see [6]).

For such a subbundle we have the decomposition
$$
(T_M \oplus T_M^*) \otimes \mathbb{C} = L \oplus \bar L, 
$$
and we may use the inner product $\langle \cdot, \cdot\rangle$ to identify $\bar L \equiv L^*$. Let 
$$
\pi_T : (T_M \oplus T_M^*) \otimes \mathbb{C} \to T_M \otimes \mathbb{C}
$$
be the projection and let $E = \pi_T (L)$. Then, the type $k = \{0, 1, \ldots, n\}$ of the generalized complex structure at $x \in M$ is defined as the codimension of $E_x \subset T_x \otimes \mathbb{C}$.

The generalized almost complex structure determined by a  symplectic structure
$$
J_\omega = 
\left(
\begin{array}{cc}
 0 & -\omega^{-1} \\
\omega & 0
\end{array}
\right)
$$
has $(+i)$-eigenbundle
$$
L = \{X - i\omega (X) \mid X \in T_M \otimes \mathbb{C}\}, 
$$
which is Courant involutive if and only if $d \omega = 0$. The type is $k = 0$ in this case. 

The generalized complex structure determined by a  complex structure 
$$
J_J = 
\left(
\begin{array}{cc}
-J & 0 \\
0 & J^\ast
\end{array}
\right)
$$
has maximal isotropic eigenbundle
$$
L = T_{0, 1} \oplus T^\ast_{1, 0},
$$
which is Courant involutive if and only if $J$ is integrable as a complex structure. The type is $k = n$ in this case (see [6]). 

Let $M$ be a smooth manifold of even dimension (i.e. $\dim_{\mathbb{R}} M = 2n$). 

A generalized complex structure $J$ is determined by its $(+i)$-eigenbundle $L \subset (T_M \oplus T_M^*)\otimes \mathbb{C}$ 
which is isotropic, satisfies $L \cap \bar L = \{0\}$ and is closed under the Courant bracket. 

To deform $J$ we will vary $L$ in the Grassmannian of maximal isotropic. Any maximal isotropic having zero intersection with $\bar L$ (this is an open set containing $L$) can be uniquely described as the graph of a homomorphism $\varepsilon : L \to \bar L$ satisfying 
$\langle \varepsilon (X), Y \rangle + \langle X, \varepsilon (Y) \rangle= 0$, $ \forall\; X, Y \in \mathcal{C}^\infty (L)$ or equivalently $\varepsilon \in \mathcal{C}^\infty(\wedge^2 L^*)$. Therefore the new isotropic is given by $L_\varepsilon=(1 + \varepsilon)L$. 
As the deformed $J$ is to remain real, we must have $\bar L_\varepsilon = (1 + \bar\varepsilon) \bar L$. Now, $L_\varepsilon$ has zero intersection with $\bar L_\varepsilon$ if and only if the endomorphism on $L \oplus L^\ast$, described by
$$
A_\varepsilon = 
\left(
\begin{array}{cc}
 1& \bar\varepsilon \\
\varepsilon & 1
\end{array}
\right)
$$
is invertible; this is the case for $\varepsilon$ in an open set around zero (see [6]). So, providing $\varepsilon$ is small enough, $J_\varepsilon = A_\varepsilon J A_\varepsilon^{-1}$ is a new generalized almost complex structure. By [11], $J_\varepsilon$ is integrable if and only if $\varepsilon \in \mathcal{C}^\infty(\wedge^2 L^\ast)$ satisfies the equation
$$
d_L \varepsilon + \frac{1}{2} [\varepsilon, \varepsilon] = 0, 
$$
which is called the Maurer-Cartan equation. Gualtieri proved in [6] a general theorem on the deformations for generalized complex structures.

\vspace{0.5cm}

\section*{4. Deformations of generalized complex structures on Kodaira surfaces.}

\vspace{0.5cm}

Let $N$ be a (primary) Kodaira surface. Recall that the tangent bundle $T_N$ is globally generated by the global vector fields
$\{ X, Y, U, V\}$ with the only non-zero Poisson bracket $[X, Y] = U$ and that the structure endomorphism $J$ is acting on $T_N$ by $JX = Y$, $JY= -X$,  $JU =V$ and $JV = -U$. Let 
$$
T = \frac{1}{2} (X - iY),\quad W = \frac{1}{2}(U - iV). 
$$

The only non-zero Poisson bracket between $T, \bar T, W, \bar W$ is
$$
[T, \bar T] = \frac{i}{2}(W + \bar W).
$$

Then, the tangent bundle $T_N$ is globally generated by $\left\{T, W, \bar T, \bar W\right\}$ and the cotangent bundle $T_N^{\ast}$ is globally generated by the dual basis of 1-forms $\left\{\omega, \rho, \bar\omega, \bar\rho\right\}$.

We have
$$
JT = iT, \; JW = iW, \; J\bar T = -i\bar T, \; J\bar W = -iW.
$$

It follows that the subbundles $T_{0,1}$ and $T_{1,0}$ of the tangent bundle $T_N$ are globally generated by $\{\bar T, \bar W\}$, respectively $\{T, W\}$ and the dual bundle $T_{1,0}^\ast$ is globally generated by $\{\omega, \rho\}$.

The standard complex structure on the Kodaira surface $N$ can be seen as a generalized complex structure given by the subbundle $L \subset (T_N \oplus T_N^{\ast})\otimes \mathbb{C}$ of the form
$$
L = \{ \bar T, \bar W, \omega, \rho\}{}^{\tilde{}} = (T_{0, 1} \oplus T_{1,0}^\ast) \otimes \mathbb{C},
$$
which is maximal isotropic and Courant involutive (see [6] or section 3).

In the following we shall study the deformations of this generalized complex structure on the Kodaira surface $N$.

We have
$$
L \cap \bar L = \{ 0 \}, \; L\otimes \bar L = (T_N \oplus T_N^{\ast})\otimes \mathbb{C}.
$$

Using the inner product
$$
\langle X+\sigma, Y+\tau \rangle = \frac{1}{2}(\sigma(Y) + \tau(X)),
$$
for $X,Y \in \mathcal{C}^{\infty}(T_N)\otimes\mathbb{C}$, $\sigma, \tau \in \mathcal{C}^{\infty}(T_N^{\ast})\otimes\mathbb{C}$, we can identify $\bar L$ with $L^{\ast}$ by the isomorphism:
$$
\theta : \bar L \displaystyle\mathop{\rightarrow}\limits^{\widetilde{}} L^{\ast},\; \theta(T) = \frac{1}{2}\omega^{\ast},\; \theta(W) = \frac{1}{2}\rho^{\ast},\; \theta(\bar\omega) = \frac{1}{2}\bar T^{\ast},\; \theta(\bar \rho) = \frac{1}{2}\bar W^{\ast}.
$$

In order to obtain the deformations of the generalized complex structure we must consider the linear maps
$$
\varepsilon : L \to \bar L,
$$
which verify the condition
$$
\langle \varepsilon(X_0), X_1\rangle + \langle X_0, \varepsilon(X_1)\rangle = 0, \; \forall X_0,X_1 \in \mathcal{C}^{\infty}(L); \leqno(*)
$$
see [6] or section 3.

Equivalently, we can consider the linear maps $\widetilde{\varepsilon} = \theta \circ \varepsilon : L \to L^{\ast}$ and, by simple computation, we have:

\vspace{0.5cm}

{\bf Lemma 4.1} {\it The map $\tilde{\varepsilon}$ is given by the following matrix:
$$
\widetilde{\varepsilon}=\frac{1}{2}
\left(
\begin{array}{cccc}
0 &  t_{32} & -t_{11} & -t_{21} \\
-t_{32} & 0 & -t_{12} & -t_{22} \\
t_{11} & t_{12} & 0 & t_{14} \\
t_{21} & t_{22} & -t_{14} & 0
\end{array}
\right),\quad t_{ij} \in \mathbb{C}.
$$}

Observe that $\widetilde{\varepsilon} \in \mathcal{C}^{\infty}(\wedge^2L^{\ast})$. From now on we shall solve the generalized Maurer-Cartan equation:
$$
d_L\widetilde{\varepsilon} + \frac{1}{2}[\tilde{\varepsilon}, \tilde{\varepsilon}] = 0,
$$
where $\tilde{\varepsilon} \in \mathcal{C}^{\infty}(\wedge^2L^{\ast})$, $[\tilde{\varepsilon}, \tilde{\varepsilon}]$ is the Schouten bracket and
$$
d_L : \mathcal{C}^{\infty}(\wedge^2L^{\ast}) \rightarrow \mathcal{C}^{\infty}(\wedge^3L^{\ast});
$$
see [6]. The derivative $d_L$ for a Lie algebroid is given in this case by the formula
$$
\begin{array}{l}
d_L\tilde{\varepsilon}(X_0,X_1,X_2) = a(X_0)\tilde{\varepsilon}(X_1,X_2) -a(X_1)\tilde{\varepsilon}(X_0,X_2) + a(X_2)\tilde{\varepsilon}(X_0, X_1) -\\[3mm]
- \tilde{\varepsilon}([X_0,X_1], X_2) + \tilde{\varepsilon}([X_0,X_2], X_1) -\tilde{\varepsilon}([X_1,X_2], X_0),
\end{array}
$$
where the anchor map $a : \mathcal{C}^{\infty}(L) \to \mathcal{C}^{\infty}(T_N)$ is the projection on the tangent bundle $T_N$ and $X_0,X_1, X_2 \in \mathcal{C}^{\infty}(L)$.

We shall use the notation:
$$
\begin{array}{l}
X_0 = u_1\bar T + u_2\bar W + u_3\omega + u_4\rho,\\[3mm]
X_1 = \alpha_1\bar T + \alpha_2\bar W + \alpha_3\omega + \alpha_4\rho,\\[3mm]
X_2 = \beta_1\bar T + \beta_2\bar W + \beta_3\omega + \beta_4\rho,
\end{array}
$$
where $u_i, \alpha_i, \beta_i \in \mathcal{C}^{\infty}(N), \; \forall i=1,2,3,4$.

We have
$$
\begin{array}{l}
\tilde{\varepsilon}(X_1,X_2) = \frac{1}{2}(\beta_1(t_{32}\alpha_2 -t_{11}\alpha_3-t_{21}\alpha_4) + \beta_2(-t_{32}\alpha_1 -t_{12}\alpha_3-t_{22}\alpha_4) +\\[3mm]
+ \beta_3(t_{11}\alpha_1+t_{12}\alpha_2+t_{14}\alpha_4) + \beta_4(t_{21}\alpha_1 +t_{22}\alpha_2 - t_{14}\alpha_3)).
\end{array}
$$

Since $a(X_0) = u_1\bar T + u_2\bar W$, by direct computation, we get

\vspace{0.5cm}

{\bf Lemma 4.2} {\it 
$$
\begin{array}{l}
a(X_0)\tilde{\varepsilon}(X_1, X_2) = \frac{1}{2}(t_{11}u_1(\bar T(\alpha_1)\beta_3 - \bar T(\beta_3)\alpha_1 - \bar T(\alpha_3)\beta_1 - \bar T(\beta_2)\alpha_3) +\\[3mm]
+ t_{12}u_1(\bar T(\alpha_2)\beta_3 + \bar T(\beta_3)\alpha_2 - \bar T(\alpha_3)\beta_2 -\bar T(\beta_2)\alpha_3) +\\[3mm]
+ t_{21}u_1(\bar T(\alpha_1)\beta_4 + \bar T(\beta_4)\alpha_1 - \bar T(\alpha_4)\beta_1 -\bar T(\beta_1)\alpha_4) +\\[3mm]
+ t_{22}u_1(\bar T(\alpha_2)\beta_4 + \bar T(\beta_4)\alpha_2 - \bar T(\alpha_4)\beta_2 - \bar T(\beta_2)\alpha_4) +\\[3mm]
+ t_{14}u_1(\bar T(\alpha_4)\beta_3 + \bar T(\beta_3)\alpha_4 - \bar T(\alpha_3)\beta_4 - \bar T(\beta_4)\alpha_3) +\\[3mm]
+ t_{32}u_1(\bar T(\alpha_2)\beta_1 + \bar T(\beta_1)\alpha_2 - \bar T(\alpha_1)\beta_2 - \bar T(\beta_2)\alpha_1) +\\[3mm]
+ t_{11}u_2(\bar W(\alpha_1)\beta_3 + \bar W(\beta_3)\alpha_1 - \bar W(\alpha_3)\beta_1 - \bar W(\beta_1)\alpha_3) +\\[3mm]
+ t_{12}u_2(\bar W(\alpha_2)\beta_3 + \bar W(\beta_3)\alpha_2 - \bar W(\alpha_3)\beta_2 - \bar W(\beta_2)\alpha_3) +\\[3mm]
+ t_{21}u_2(\bar W(\alpha_1)\beta_4 + \bar W(\beta_4)\alpha_1 - \bar W(\alpha_4)\beta_1 - \bar W(\beta_1)\alpha_4) +\\[3mm]
+ t_{22}u_2(\bar W(\alpha_2)\beta_4 + \bar W(\beta_4)\alpha_2 - \bar W(\alpha_4)\beta_2 - \bar W(\beta_2)\alpha_4) +\\[3mm]
+ t_{14}u_2(\bar W(\alpha_4)\beta_3 + \bar W(\beta_3)\alpha_4 - \bar W(\alpha_3)\beta_4 - \bar W(\beta_4)\alpha_3) +\\[3mm]
+ t_{32}u_2(\bar W(\alpha_2)\beta_1 + \bar W(\beta_1)\alpha_2 - \bar W(\alpha_1)\beta_2 - \bar W(\beta_2)\alpha_1)).
\end{array}
$$}

Analogous formulae can be written for the terms $a(X_1)\tilde{\varepsilon}(X_0,X_2)$ and $a(X_2)\tilde{\varepsilon}(X_0,X_1)$.

We shall use the following notation:
$$
X_0 = X + \sigma, \; X_1 = Y + \tau, \; X_2 = Z + \delta,
$$
where
$$
X = u_1\bar T + u_2\bar W,\; Y = \alpha_1\bar T + \alpha_2\bar W,\; Z = \beta_1\bar T + \beta_2\bar W \in \mathcal{C}^{\infty}(T_N)
$$
and
$$
\sigma = u_3\omega + u_4\rho,\; \tau = \alpha_3\omega + \alpha_4\rho,\; \delta = \beta_3\omega + \beta_4\rho \in \mathcal{C}^{\infty}(T_N^{\ast}).
$$

It is easy to see that:
$$
\mathcal{L}_X\tau =( X(\alpha_3) - \frac{i}{2}u_1\alpha_4)\omega + X(\alpha_4)\rho,
$$
$$
\mathcal{L}_Y\sigma = (Y(\alpha_3) - \frac{i}{2}\alpha_1 u_4)\omega + Y(u_4)\rho,
$$
and
$$
i_X\tau = 0, \; i_Y\sigma = 0.
$$

Using the definition of the Courant bracket, by direct computation, we get:

\vspace{0.5cm}

{\bf Lemma 4.3} {\it  
$$
\begin{array}{l}
[X_0, X_1] = (X(\alpha_1) - Y(u_1))\bar T + (X(\alpha_2) - Y(u_2))\bar W +\\[3mm]
+ (X(\alpha_3) - Y(u_3) + \frac{i}{2}(\alpha_1u_4 - u_1\alpha_4))\omega + (X(\alpha_4) - Y(u_4))\rho.
\end{array}
$$}

We have analogous formulae for $[X_0, X_2]$ and $[X_1, X_2]$.

Now, a tedious but direct computation gives:

\vspace{0.5cm}

{\bf Lemma 4.4} {\it  
$$
\begin{array}{l}
\tilde{\varepsilon}([X_0,X_2],X_1)= \frac{1}{2}(t_{11}(-\alpha_1u_1\bar T(\beta_3) - \alpha_1u_2\bar W(\beta_3) + \alpha_1\beta_1\bar T(u_3) +\\[3mm] 
+ \alpha_1\beta_2\bar W(u_3) - \frac{i}{2}\alpha_1\beta_1u_4 + \frac{i}{2}\alpha_1u_1\beta_4 + \alpha_3u_1\bar T(\beta_1)+ \\[3mm]
+ \alpha_3u_2\bar W(\beta_1) - \alpha_3\beta_1\bar T(u_1) - \alpha_3\beta_2\bar W(u_1)) + t_{12}(-\alpha_2u_1\bar T(\beta_3) -\\[3mm]
- \alpha_2u_2\bar W(\beta_3) + \alpha_2\beta_1\bar T(u_3) + \alpha_2\beta_2\bar W(u_3) - \frac{i}{2}\alpha_2\beta_1u_4 +\\[3mm]
+ \frac{i}{2}\alpha_2u_1\beta_4 + \alpha_3u_1\bar T(\beta_2) + \alpha_3u_2\bar W(\beta_2) - \alpha_3\beta_1\bar T(u_2) -\\[3mm]
- \alpha_3\beta_2\bar W(u_2)) + t_{21}(-\alpha_1u_1\bar T(\beta_4) - \alpha_1u_2\bar W(\beta_4) + \alpha_1\beta_1\bar T(u_4) +\\[3mm]
+ \alpha_1\beta_2\bar W(u_4) + \alpha_4u_1\bar T(\beta_1) + \alpha_4u_2\bar W(\beta_1) - \alpha_4\beta_1\bar T(u_1) -\\[3mm]
- \alpha_4\beta_2\bar W(u_1)) + t_{22}(-\alpha_2u_1\bar T(\beta_4) - \alpha_2u_2\bar W(\beta_4) + \alpha_2\beta_1\bar T(u_4) +\\[3mm]
+ \alpha_2\beta_2\bar W(u_4) + \alpha_4u_1\bar T(\beta_2) + \alpha_4u_2\bar W(\beta_2) - \alpha_4\beta_1\bar T(u_2) -\\[3mm]
- \alpha_4\beta_2\bar W(u_2)) + t_{14}(\alpha_3u_1\bar T(\beta_4) + \alpha_3u_2\bar W(\beta_4) - \alpha_3\beta_1\bar T(u_4) -\\[3mm]
- \alpha_3\beta_2\bar W(u_4) - \alpha_4u_1\bar T(\beta_3) - \alpha_4u_2\bar W(\beta_3) + \alpha_4\beta_1\bar T(u_3) +\\[3mm]
+ \alpha_4\beta_2\bar W(u_3) - \frac{i}{2}\alpha_4\beta_1u_4 + \frac{i}{2}\alpha_4u_1\beta_4) + t_{32}(\alpha_1u_1\bar T(\beta_2) +\\[3mm]
+ \alpha_1u_2\bar W(\beta_2) - \alpha_1\beta_1\bar T(u_2) - \alpha_1\beta_2\bar W(u_2) - \alpha_2u_1\bar T(\beta_1) -\\[3mm]
- \alpha_2u_2\bar W(\beta_1) + \alpha_2\beta_1\bar T(u_1) + \alpha_2\beta_2\bar W(u_1))).
\end{array} 
$$}

There are analogous formulae for the terms $\tilde{\varepsilon}([X_0, X_1], X_2)$ and $\tilde\varepsilon([X_1, X_2], X_0)$.

By all the computations above, we get the first result:

\vspace{0.5cm}

{\bf Theorem 4.5} {\it The differential $d_L$ of $\tilde\varepsilon$ is given by the following formula:
$$
(d_L \tilde\varepsilon)(X_0, X_1, X_2) = 
\frac{i}{2} t_{12} (\alpha_4\beta_1u_2 - \alpha_1\beta_4u_2 - 
\alpha_2\beta_1u_4 + \alpha_2\beta_4u_1 + \alpha_1\beta_2u_4 - 
\alpha_4\beta_2u_1),
$$
where $X_0, X_1, X_2 \in \mathcal{C}^\infty (L)$.}

Now, we have to compute the Schouten bracket $[\tilde\varepsilon, \tilde\varepsilon]$. We can write the omomorphism $\tilde\varepsilon$ in the following form:
$$
\tilde\varepsilon =t_{32} \bar T^\ast \wedge \bar W^\ast - t_{11}\bar T^\ast \wedge \omega^\ast - t_{21} \bar T^\ast 
\wedge \rho^\ast - t_{12}\bar W^\ast \wedge \omega^\ast - t_{22} \bar W^\ast \wedge \rho^\ast + t_{14} \omega^\ast \wedge \rho^\ast.
$$

The next result is the following:

\vspace{0.5cm}

{\bf Theorem 4.6} {\it For the Schouten bracket we get:}
$$
[\tilde\varepsilon, \tilde\varepsilon] = 0.
$$

\vspace{0.2cm}

{\it Proof.} It suffices to prove that the Schouten brackets of all pairs of the set 
$\{ \bar T^\ast \wedge \bar W^\ast, \bar T^\ast \wedge \omega^\ast, \bar T^\ast \wedge \rho^\ast, \bar W^\ast \wedge \omega^\ast, 
\bar W^\ast \wedge \rho^\ast, \omega^\ast \wedge \rho^\ast\}$ are zero. 

By using the identifications
$$
\bar T^\ast \equiv  2 \bar\omega, \quad \bar W^\ast \equiv 2\bar\rho, \quad \omega^\ast \equiv 2 T, \quad \rho^\ast \equiv 2W
$$
and the definitions of the Schouten and Courant brackets we get the result after tedious computations. We shall present only some typical cases. The first case is the following:
$$
\begin{array}{l}
 [\bar T^\ast \wedge \bar W^\ast, \bar T^\ast \wedge \bar W^\ast] =  16[ \bar\omega \wedge \bar\rho, \bar \omega \wedge \bar\rho]=\\
= 16([\bar\omega, \bar\omega]\wedge \rho\wedge\rho - 
[\rho, \omega]\wedge \omega\wedge \rho - 
[\omega, \rho] \wedge\rho\wedge\omega +
[\rho, \rho] \wedge \omega \wedge\omega] = 0,
\end{array}
$$
since
$$
\rho \wedge \rho = 0, \quad \omega \wedge \omega =0, \quad [\omega, \rho] = 0.
$$
The second case:
$$
[\bar T^\ast \wedge \bar W^\ast,  \bar T^\ast \wedge \omega^\ast] = 16[\bar\omega \wedge \bar \rho, \bar\omega \wedge T] =
-16 [\bar\omega, T] \wedge \bar\rho \wedge\bar\omega.
$$

Now, $[\bar\omega, T] = - \mathcal{L}_T \bar\omega$ and for any 
$$
\widetilde Y = \gamma_1 \bar T + \gamma_2 \bar W + \gamma_3 T + \gamma_4 W \in \mathcal{C}^\infty(T_N), 
$$
we get
$$
(\mathcal{L}_T \bar\omega)(\tilde Y) = T(\gamma_1) - \bar\omega(\mathcal{L}_T\tilde Y).
$$

But
$$
\mathcal{L}_T \widetilde Y = -\frac{i}{2} \gamma_1 (\bar W + W) + T(\gamma_1)\bar T + T(\gamma_2) \bar W +T(\gamma_3) T + T(\gamma_n)W,
$$
hence 
$$
(\mathcal{L}_T \bar\omega)(\widetilde Y) = T(\gamma_1) - T(\gamma_1) = 0.
$$

It follows $[\bar\omega, T] = 0$.

The third case:
$$
\begin{array}{l}
 [\bar T^\ast \wedge \omega^\ast, \bar W^\ast \wedge \rho^\ast] = 16 [\bar\omega \wedge T, \bar \rho\wedge W]=\\
16([\bar\omega, \bar\rho]\wedge T\wedge W -[\bar\omega, W]\wedge T\wedge \bar\rho - 
[T, \bar\rho]\wedge \bar\omega \wedge W + [T, W] \wedge \bar\omega \wedge \bar\rho).
\end{array}
$$

Since
$$
[\bar\omega, \bar\rho] = 0, \quad [\bar\omega, W] = 0, \quad [T, W] = 0, \quad [T, \bar\rho] = -\frac{i}{2}\bar\omega, 
$$
we get
$$
[\bar T^\ast \wedge \omega^\ast, \bar W^\ast \wedge \rho^\ast] = +16 \frac{i}{2} \bar\omega \wedge \bar\omega \wedge W = 0.
$$

From the two theorems, we get the following:

\vspace{0.5cm}

{\bf Corollary 4.7 } {\it The solutions of the generalized Maurer-Cartan equation are given by}
$$
\bar\varepsilon = 4(t_{32} \bar\omega \wedge \bar\rho - t_{11} \bar\omega \wedge T - 
t_{21} \bar\omega \wedge W - t_{22} \bar\rho \wedge W + t_{14} T \wedge W)
$$

Now, we need the following result: 

\vspace{0.5cm}

{\bf Proposition 4.8} {\it The image of the differential
$$
d_L : \mathcal{C}^\infty (L^\ast) \to \mathcal{C}^\infty (\wedge^2 L^*)
$$
is generated by $\bar\omega \wedge W$. }

\vspace{0.5cm}

{\it Proof.} Let $\sigma = t_1 \bar T^\ast+ t_2 \bar W^\ast + t_3 \omega^\ast + t_4 \rho^\ast \in \mathcal{C}^\infty (L^\ast)$ and 
$ X_0 = u_1 \bar T + u_2 \bar W + u_3 \omega + u_4 \rho$, $X_1 = \alpha_1 \bar T + \alpha_2 \bar W + \alpha_3 \omega + \alpha_4 \rho \in \mathcal{C}^\infty (L)$, where $t_i, u_i, \alpha_i \in \mathcal{C}^\infty(N)$, $i = 1,2,3,4$. We have:
$$
(d_L \sigma)(X_0, X_1) = a(X_0) \sigma(X_1) - a(X_1) \sigma (X_0) - 
\sigma([X_0, X_1]).
$$

Similar, but shorter, computation as in the proof of Theorem 4.5 gives us the formula:
$$
(d_L \sigma)(X_0, X_1) = \frac{i}{2} t_3 (\alpha_1u_4 - u_1\alpha_4) = i t_3 \bar T^\ast \wedge \rho^\ast(X_0, X_1).
$$

It follows that the image of the differential $d_L : \mathcal{C}^\infty (L^\ast) \to \mathcal{C}^\infty (\wedge^2 L^\ast)$ is generated by 
$\bar T^\ast \wedge\rho^\ast$ or, equivalently, by $\bar \omega \wedge W$.

Since deformations of generalized complex structures, which differ by an element in the image of the  differential $d_L : \mathcal{C}^\infty (L^\ast) \to \mathcal{C}^\infty (\wedge^2 L^\ast)$, are equivalent, we get the main result of the paper:

\vspace{0.5cm}

{\bf Theorem 4.9} {\it The deformations in the sense of generalized complex structures of the standard complex structure on a (primary) Kodaira surface $N$ are given by 
$$
\tilde\varepsilon = t_{32}\bar\omega \wedge \bar\rho - t_{11} \bar\omega \wedge T - 
t_{22} \bar\rho \wedge W + t_{14} T\wedge W,
$$
where $(t_{32}, t_{11}, t_{22}, t_{14})\in \mathbb{C}^4$. }

\vspace{0.5cm}

{\bf Remark} The above result shows that the deformations in the sense of generalized complex structures of the
standard complex structure on a (primary) Kodaira surface are the same as the extended deformations (in the sense of Kontsevich [9], Barannikov-Kontsevich [1]) in degree two, obtained by Poon [12] using differential Gerstenhaber algebras. In particular, taking the parameters $t_{32} = 0$ and $t_{14} = 0$ we get the classical deformations of complex structures obtained by Borcea [2].

We have the following:

\vspace{0.3cm}

{\bf Corollary 4.10} {\it The family of deformations of generalized complex structures on a (primary) Kodaira surface $N$, 
given by 
$$
\widetilde\varepsilon = t_{32} \bar\omega \wedge \bar\rho - t_{11}\bar\omega \wedge T - t_{22} \bar\rho \wedge W + t_{14} T \wedge W, 
$$
with $(t_{32}, t_{11}, t_{22}, t_{14}) \in U \subset \mathbb{C}^4$, where $U$ is an  open neighborhood of $0 \in \mathbb{C}^4$, is a smooth locally complete family. }

\vspace{0.5cm}

{\it Proof.} By Theorem 4.6 we have $[\widetilde\varepsilon, \widetilde\varepsilon] = 0$. From the  definition of the obstruction 
map $\phi$  given in the proof of the Theorem 5.4 in [6] (see also [10]) we get $\phi = 0$. Then, applying again Theorem 5.4 in [6] it follows that the above family of deformations is a smooth locally complete family in a open neighborhood $U$ of $0 \in \mathbb{C}^4$.

The generalized complex structure defined by the  deformation given by the homomorphism
$$
\varepsilon : L \to \widetilde L \simeq L^\ast, \quad (\varepsilon \in \mathcal{C}^\infty (\wedge^2 L^\ast))
$$
is, uniquely, described by the new  isotropic subbundle
$$
L_\varepsilon  = (1 +\varepsilon) L \subset (T_N \oplus T_N^\ast) \otimes \mathbb{C}.
$$

We have the following result:

\vspace{0.5cm}

{\bf Proposition 4.11} {\it Let $N$ be a (primary) Kodaira surface. The type of the generalized complex structure given by the subbundle $L_\varepsilon \subset (T_N \oplus T_N^\ast) \otimes \mathbb{C}$ is $k_\varepsilon = 0$ (symplectic type) or $k_\varepsilon = 2$ (complex type). }

\vspace{0.5cm}

{\it Proof.} The type $k_\varepsilon$ of a generalized complex structure is the codimension in any fibre of the projection of $L_\varepsilon$ on 
$T_N \otimes \mathbb{C}$ by the canonical map
$$
p_1 : (T_N \oplus T_N^\ast) \otimes \mathbb{C} \to T_N \otimes \mathbb{C}.
$$

Since
$$
L= \{ \bar T, \bar W, \omega, \rho\}^{\tilde{}} = (T_{0,1} \oplus T_{1, 0}^\ast) \otimes \mathbb{C}, 
$$
we get:
$$
\begin{array}{l}
 (1 + \varepsilon)(\bar T)= \bar T + t_{11} T- t_{32} \bar\rho\\
(1 + \varepsilon) (\bar W)= \bar W + t_{22} W - t_{32} \bar\omega\\
(1 + \varepsilon) (\omega) = - t_{14} W+ \omega - t_{11} \bar\omega\\
(1 + \varepsilon) (\rho) = t_{14} T + \rho - t_{22} \bar\rho.
\end{array}
$$

It follows that the projection of $L_\varepsilon$ on $T_N \otimes \mathbb{C}$ is globally generated by 
$$
\{\bar T + t_{11}T, \bar W + t_{22} W, -t_{14}W, t_{14}T\}.
$$

If $t_{14} \neq 0$, then the type $k_\varepsilon = 0$ (symplectic type) and, if $t_{14} = 0$,then $k_\varepsilon = 2$ (complex type). 

\vspace{0.5cm}

{\bf Remark} If, in the case $t_{14} = 0$, we have also $t_{32} = 0$, we get classical deformations of complex structures. If,  in the case $t_{14} = 0$, we have $t_{32} \neq 0$, we get examples of generalized complex structures of complex type, which are not classical complex structures.


\vspace{1cm}

\begin{thebibliography}{ll}
\bibitem{[1]} S. Barannikov, M. Kontsevich, Frobenius manifolds and formality of Lie algebras of polyvector fields, Internat. Math. Res. Notes 14(1998), 201-215.

\bibitem{[2]} C. Borcea, Moduli for Kodaira surfaces, Comp. Math. 52(1984), 373-380.

\bibitem{[3]} V. Br\^\i nz\u anescu, Holomorphic Vector Bundles over Compact Complex Surfaces, Lecture Notes in Math. 1624, Springer 1996. 

\bibitem{[4]} T. Courant, Dirac manifolds, Trans. Amer. Math. Soc. 319(1990), 631-661.

\bibitem{[5]} T. Courant, A. Weinstein, Beyond Poisson structures, In Action hamiltoniennes de groupes, Troisi\` eme th\' eor\` eme de Lie (Lyon,1986), volume 27 of Travaux en Cours, 39-49, Hermann,Paris 1988.

\bibitem{[6]} M.Gaultieri, Generalized complex geometry, D. Phil. Thesis, St. John's College, University of Oxford, 2003. 

\bibitem{[7]} N. Hitchin, Generalized Calabi-Yau manifolds, Q. J. Math. 54 (3) (2003), 281-308. 

\bibitem{[8]} K. Kodaira, On the structure of compact complex analytic surfaces I, Amer. J. Math. 86(1964), 751-798.

\bibitem{[9]} M. Kontsevich, Homological algebra of mirror symmetry, In Proceedings of the International Congress of Mathematicians, Vol. 1, 2(Z\" urich, 1994), pages 120-139, Basel, 1995, Birkh\" auser.

\bibitem{[10]} M. Kuranishi, New proof for the existence of locally complete families of complex structures, In Proc. Conf. Complex Analysis (Mineapolis, 1964), 142-154, Springer 1965. 

\bibitem{[11]} Z.-J. Liu, A. Weinstein, P. Xu, Manin triples for Lie bialgebroids, J. Diff. Geom. 45,(1997),  547-574.

\bibitem{[12]} Y.S. Poon, Extended deformation of Kodaira surfaces, J. reine angew. Math. 590(2006), 45-65.



\end{thebibliography}
\end{document}